\documentclass[11pt]{article}
\usepackage{graphicx}
\usepackage{amssymb}
\evensidemargin \oddsidemargin
\def\x{\mbox{x}}

\def\b#1{{\mathbb #1}}
\def\c#1{{\cal #1}}
\def\nn{\nonumber \\}
\def\1{{\bf 1}}
\def\RH{\mbox{$\hat {\sf R}$\,}}
\def\PH{\mbox{$\sf P$}}
\def\P{\mbox{$\cal P$}}

\def\A{\mbox{$\cal A$}}
\def\id{\mbox{id\,}}

\newcommand{\be}{\begin{equation}}
\newcommand{\ee}{\end{equation}}
\newcommand{\bea}{\begin{eqnarray}}
\newcommand{\eea}{\end{eqnarray}}
\newcommand{\ba}{\begin{array}}
\newcommand{\ea}{\end{array}}

\begin{document}
\title{$q$-Deformed quaternions and $su(2)$ 
instantons\footnote{Talk given at ``Noncommutative Geometry in 
Field and String Theories'',
Satellite Workshop of ``Corfu Summer Institute 2005'',
18-20/9/05,  Corfu, Grecia. Appeared in the Proceedings. }}

\author{        Gaetano Fiore,  \\ ~ \\
        Dip. di Matematica e Applicazioni, Universit\`a ``Federico II''\\
        V. Claudio 21, 80125 Napoli, Italy\\
        \and
        \and
        I.N.F.N., Sezione di Napoli,\\
        Complesso MSA, V. Cintia, 80126 Napoli, Italy
        }
\date{}

\maketitle
\abstract{
We have recently introduced the notion of a $q$-quaternion bialgebra and
shown its strict link with the $SO_q(4)$-covariant quantum Euclidean
space $\b{R}_q^4$. Adopting the available differential geometric tools on the
latter and the quaternion language we have formulated and found 
solutions of the (anti)selfduality equation [instantons and
multi-instantons] of a would-be deformed $su(2)$ Yang-Mills theory on this
quantum space. The solutions depend on some noncommuting parameters,
indicating that the moduli space of a complete theory should be a
noncommutative manifold. We summarize these results and add an explicit
comparison between the two $SO_q(4)$-covariant differential calculi on
$\b{R}_q^4$ and the two 4-dimensional bicovariant differential calculi on the
bi- (resp. Hopf) algebras $M_q(2),GL_q(2),SU_q(2)$,
showing that they essentially coincide.}


\section{Introduction}
The construction of gauge field theories on noncommutative manifolds
has been  the subject of quite a lot of work in recent years. A crucial test
of it is the search of instantonic solutions, especially
after the discovery \cite{NekSch98}  that deforming $\b{R}^4$ into
the Moyal-Weyl noncommutative Euclidean space $\b{R}_{\theta}^4$ regularizes
the zero-size singularities of the instanton moduli space (see also
\cite{SeiWit99}). 
Various other noncommutative geometries have been considered 
(see e.g. \cite{ConLan01,BonCicTar00,DabLanMas01,LanvSu06}). They do not
always completely fit Connes' standard framework of noncommutative geometry
\cite{Con94}, thus stimulating attempts of generalizations. Among the available
deformations of $\b{R}^4$ there is also the Faddeev-Reshetikhin-Takhtadjan
noncommutative Euclidean space $\b{R}_q^4$ covariant under $SO_q(4)$
\cite{FadResTak89}. This, as other quantum group covariant
noncommutative spaces (shortly: quantum spaces), is maybe even more
problematic for the formulation \cite{JurMoeSchSchWes01} of a gauge field
theory on  like $\b{R}_q^4$. One main reason is the lack of a proper (i.e.
cyclic) trace to define gauge invariant observables (action, etc). Another one
is the $\star$-structure of the differential calculus, which for real $q$ is
problematic. Nevertheless, in our main Ref. \cite{Fio05i} we have left these
two issues aside and investigated about (anti)selfduality equations on it
and their solutions. Here we summarize these results adding some detail.

As a first step 
we recall our notion  \cite{Fio05i} of a $q$-deformed
quaternion as the defining matrix of a copy of $SU_q(2)\times
\b{R}^{\ge}$ ($\b{R}^{\ge}$ denoting the semigroup of nonnegative
real numbers), or equivalently of the $2\times 2$ defining quantum matrix 
of $M_q(2)$ endowed with the same $\star$-structure of $SU_q(2)$ (more details
will be given in \cite{Fio05q}), and  that
its entries can be regarded also as coordinates of $\b{R}_q^4$.  As
on ordinary $\b{R}^4$, this will much simplify the search and classification of
instantons in Yang-Mills theory. We also recall that the quantum sphere $S_q^4$ 
of \cite{DabLanMas01}
can be regarded as a compactification of the corresponding $\star$-algebra.
We then show that the two $SO_q(4)$-covariant
differential calculi on $\b{R}_q^4$ \cite{CarSchWat91} 
coincide with the two 4-dimensional bicovariant differential calculi 
\cite{S2,SchWatZum92} on the
bi- (resp. Hopf) algebras $M_q(2),GL_q(2)$, so that upon imposing 
the unit $q$-determinant condition one obtains Woronowicz pioneering $4D\pm$  
bicovariant differential calculi \cite{Wor89,PodWor90} on $SU_q(2)$
(this had been only announced in \cite{Fio05i}).
Using the
Hodge duality map \cite{Fio94,Fio04JPA} on $\b{R}_q^4$ in $q$-quaternion
language we have formulated (anti)self-duality equations and found
\cite{Fio05i} solutions $A$, in the form of 1-form valued $2\times 2$
matrices, that closely resemble their undeformed counterparts (instantons) in
$su(2)$ Yang-Mills theory on $\b{R}^4$. [The (still missing)
complete gauge theory might be however  a deformed $u(2)$ rather
than $su(2)$ Yang-Mills theory.]. The projector characterizing the instanton
projective module (playing the role of the vector bundle) of \cite{DabLanMas01}
in $q$-quaternion language takes exactly the same natural form as in
the undeformed theory.
The ``coordinates of the center'' of the instanton are nevertheless
noncommuting parameters, differently from the Nekrasov-Schwarz theory. We have
also found multi-instantons solutions: they are again parametrized by
noncommuting parameters playing the role of ``size'' and
``coordinates of the center'' of the (anti)instantons. This
indicates that the moduli space of a complete theory should be a
noncommutative manifold. This is  similar to what was proposed
in \cite{IvaLecMue04} for $\b{R}_{\theta}^4$ for selfdual
deformation parameters $\theta_{\mu\nu}$.

\section{The $q$-quaternion bialgebra $C(\b{H}_q)$}

Any element $X$ in the (undeformed) quaternion algebra $\b{H}$ is given
by
\[ X=\x_1+\x_2i+\x_3j+\x_4k,
\] with $\x\in \b{R}^4$ and imaginary $i,j,k$ fulfilling
\[i^2=j^2=k^2=-1,\quad ijk=-1.
\]
Replacing $i,j,k$ by Pauli matrices $\times$ imaginary unit ${\rm
i}$ we get
\[ X\leftrightarrow x\equiv\left(\ba{cc}
\x_1+\x_4{\rm i} \:& \x_3+\x_2{\rm i}\\
-\x_3+\x_2{\rm i} \:& \x_1-\x_4{\rm i} \ea\right)=: \left(\ba{cc}
\alpha\:& \gamma\\
-\gamma^{\star}\:& \alpha^{\star} \ea\right) \] (where $\alpha,
\gamma\in\b{C}$), and the quaternionic product becomes represented
by matrix multiplication. Therefore  $\b{H}$ essentially consists
 of all complex $2\times 2$ matrices of this
form.

This can be $q$-deformed as follows.
We just pick the pioneering definition of the (Hopf) $*$-algebra
$C\left(SU_q(2)\right)$ \cite{Wor87,Wor87a} without imposing the det$_q$=1
condition: for $q\in\b{R}$ consider the unital associative
$\star$-algebra $\A\equiv C(\b{H}_q)$ generated by elements $\alpha,
\gamma,\alpha^{\star},\gamma^{\star}$ fulfilling the commutation
relations \be \ba{l}
\alpha\gamma=q\gamma\alpha,\qquad\alpha\gamma^{\star}=q\gamma^{\star}\alpha,
\qquad\gamma\alpha^{\star}=q\alpha^{\star}\gamma,\\[8pt]
\gamma^{\star}\alpha^{\star}=q\alpha^{\star}\gamma^{\star},\qquad
[\alpha,\alpha^{\star}]=(1\!-\!q^2)\gamma\gamma^{\star} \qquad
[\gamma^{\star},\gamma]=0.\ea \label{explqquatcomrel} \ee
Introducing the matrix
\[
x\equiv \left(\ba{ll} x^{11}\: & x^{12}\\
x^{21}\:& x^{22}\ea\right):=\left(\ba{cc} \alpha\:& -q\gamma^{\star}\\
\gamma\: & \alpha^{\star} \ea\right)
\]
 we can rewrite these commutation relations as 
\be 
\hat R x_1x_2=x_1x_2\hat R
      \label{qquatcomrel}
\ee 
and the conjugation relations as
$x^{\alpha\beta}{}^{\star}=\epsilon^{\beta\gamma}x^{\delta\gamma}
\epsilon_{\delta\alpha}$, i.e. 
\be 
x^{\dagger}=\bar x\qquad\qquad
\mbox{where } \bar a:=\epsilon^{-1} a^T\epsilon \quad\forall a\in
M_2. \label{qquatstarrel} 
\ee 
Here we have used the $\epsilon$-tensor and
the braid matrix of $M_q(2),GL_q(2),SU_q(2)$,
\be
\epsilon\!=\!\left(\ba{cc} 0\: & 1\\ -q \: & 0 \ea\right)\!=\!-q
\epsilon^{-1},\qquad\qquad \hat R^{\alpha\beta}_{\gamma\delta}
=q\delta^{\alpha}_{\gamma}\delta^{\beta}_{\delta}+\epsilon^{\alpha\beta}
\epsilon_{\gamma\delta} .     \label{qepsilon}
\ee
[with $\epsilon\!\equiv\!(\epsilon_{\alpha\beta})$ and
$\epsilon^{-1}\!\equiv\!(\epsilon^{\alpha\beta})$]; note that 
$\hat R^T=\hat R$. So
$\A:=C(\b{H}_q)$ can be endowed also with a bialgebra  structure (we
are not excluding the possibility that $x\equiv {\bf 0}_2$), more
precisely a real section of the bialgebra $C\left(M_q(2)\right)$ of
$2\times 2$ quantum matrices \cite{Dri86,Wor87a,FadResTak89}. Since the coproduct
$$
\Delta(x^{\alpha\gamma})=(ax)^{\alpha\gamma}
$$
is an algebra map,  the matrix product $ax$ of any two matrices
$a,x$ with mutually commuting entries and fulfilling
(\ref{qquatcomrel}-\ref{qquatstarrel}) again fulfills the latter.
Therefore we shall call any such matrix $x$ a {\it $q$-quaternion},
and $\A:=C(\b{H}_q)$ the $q$-quaternion bialgebra.

As well-known, the socalled `$q$-determinant' of $x$ \be |x|^2
\equiv\det{}_q(x):=x^{11}x^{22}-qx^{12}x^{21}=\alpha^{\star}\alpha +
\gamma^{\star}\gamma\sim x^{\alpha\alpha'}x^{\beta\beta'}
\epsilon_{\alpha\beta}\epsilon_{\alpha'\beta'}, \ee is central,
manifestly nonnegative-definite and group-like. It is zero iff
$x$ is. Relations (\ref{qquatcomrel}) can be also equivalently
reformulated as
\be
x\bar x=\bar xx=|x|^2I_2                \label{blu}
\ee
($I_2$ denotes the unit $2\times 2$ matrix).
If we extend $C(\b{H}_q)$ assuming the existence of a new (central,
positive-definite) generator $|x|^{-1}$ (this will imply that $x$ cannot
vanish at the representation level),  one finds that $x$ is invertible with
inverse 
\be
x^{-1}=\frac{\bar x}{|x|^2}.
\ee
$C(\b{H}_q)$ becomes a Hopf
$\star$-algebra  [a real section of $C\left(GL_q(2)\right)$].
The matrix elements
of $T:=\frac x{|x|}$ fulfill the relations (\ref{qquatcomrel}) and
\be
T^{\dagger}=T^{-1}=\overline{T}, \qquad\qquad\det{}_q(T)=\1, \label{starSUq2}
\ee
namely generate as a quotient algebra $C\left(SU_q(2)\right)$
\cite{Wor87,Wor87a}, therefore in this case the entries of $x$
generate  the
(Hopf) $\star$-algebra of functions on the quantum group $SU_q(2)\times
GL^+(1)$, in analogy with the $q=1$ case.

\section{Identification of $\b{H}_q$ with $\b{R}^4_q$, and links with other
algebras}

As a $\star$-algebra, $\A:=C(\b{H}_q)$ coincides with the algebra of
functions on the $SO_q(4)$-covariant quantum Euclidean Space
$\b{R}_q^4$ of \cite{FadResTak89}, identifying their generators as
\be
x^1= qx^{11},\quad x^2= x^{12},\quad x^3=- qx^{21},\quad x^4= x^{22}.
                                                    \label{lintra}
\ee
We shall denote by
$B\equiv(B^a_{\alpha\alpha'})$  this (diagonal and invertible) matrix entering
the linear transformation $x^a=B^a_{\alpha\alpha'}x^{\alpha\alpha'}$.
We illustrate the relation between the two starting 
from the braid matrix of $SO_q(4)$, which is obtained as 
\be
\RH\equiv\big(\RH^{ab}_{cd}\big)=q^{-1}\c{B}\big(\hat R\otimes_{\b{C}}
\hat R\big)\c{B}^{-1}
\ee
($\RH$ fulfills the braid equation because $\hat R$ does),
where $\c{B}^{ab}_{\alpha\beta\alpha'\beta'}:=B^a_{\alpha\alpha'}
B^b_{\beta\beta'}$. Its decomposition
\be
\RH = q{\sf P\!}_s - q^{-1}{\sf P\!}_A + q^{-3}{\sf P\!}_t       \label{projectorR}
\ee
in orthogonal projectors follows from that of the braid matrix of
$GL_q(2)$,
\be
\hat R=q\P_s-q^{-1}\P_a,
\ee
since ${\sf P\!}:=\c{B}(\P\otimes_{\b{C}}\c{P}')\c{B}^{-1}$
is a projector whenever $\P,\c{P}'$ are\footnote{The orthonormality relations for
the $\c{P}_{\mu}$, with $\mu=s,a$,
\be
\c{P}_{\mu}\c{P}_{\nu} = \c{P}_\mu \delta_{\mu\nu}, \qquad
\sum_{\mu}\c{P}_{\mu}= I,
\ee
trivially imply the orthogonality relations for the
${\sf P\!}_{\mu}$, with $\mu=s,a,a',t$.}. In fact,  
\be
\ba{ll}
{\sf P\!}_s=\c{B}(\P_s\!\otimes_{\b{C}}\!\P_s) \c{B}^{-1},
\qquad\qquad &{\sf P\!}_t=\c{B}(\P_a\!\otimes_{\b{C}}\! \P_a)\c{B}^{-1},\\[8pt]
{\sf P\!}_a=\c{B}(\P_s\!\otimes_{\b{C}}\! \P_a)\c{B}^{-1},
\qquad\qquad &{\sf P\!}_{a'}=\c{B}(\P_a\!\otimes_{\b{C}}\! \P_s)\c{B}^{-1},\\[8pt]
{\sf P\!}_A={\sf P\!}_a+{\sf P\!}_{a'}.
\ea                                \label{defproj}
\ee
$\c{P}_s$, $\c{P}_a$, are respectively $GL_q(2)$-covariant
deformations of the symmetric and
antisymmetric projectors, and have dimension 3,1.
They can be expressed in terms
of the $q$-deformed $\epsilon$-tensor by
\be
\c{P}_a{}^{\alpha\beta}_{\gamma\delta}=-
\frac{\epsilon^{\alpha\beta}\epsilon_{\gamma\delta}}{q+q^{-1}},
\qquad\qquad
\c{P}_s{}^{\alpha\beta}_{\gamma\delta}=\delta^{\alpha}_{\gamma}
\delta_{\delta}^{\beta}+
\frac{\epsilon^{\alpha\beta}\epsilon_{\gamma\delta}}{q+q^{-1}}.
\label{exproj}
\ee
${\sf P\!}_s$, ${\sf P\!}_A$, ${\sf P\!}_t$
are $SO_q(4)$-covariant deformations of the symmetric
trace-free, antisymmetric and trace projectors respectively;
as we shall see ${\sf P\!}_a,{\sf P\!}_{a'}$ are projectors respectively on the
selfdual and antiselfdual 2-forms subspaces. By (\ref{defproj})
${\sf P\!}_s,{\sf P\!}_a,{\sf P\!}_{a'},{\sf P\!}_A,{\sf P\!}_t$ respectively have dimensions
9,3,3,6,1, and
\be
{\sf P\!}_t{}_{kl}^{ij} = (g^{sm}g_{sm})^{-1} g^{ij}g_{kl}
= \frac1{(q+q^{-1})^2} g^{ij}g_{kl}
\label{Pt}
\ee
where the $4 \!\times\! 4$ matrix $g_{ab}$ (denoted as $C_{ab}$ in
\cite{FadResTak89}) is given by
\be
g_{ab}=B^{-1}{}^{\alpha\alpha'}_aB^{-1}{}^{\beta\beta'}_b
\epsilon_{\alpha\beta}\epsilon_{\alpha'\beta'};
\ee
it is the $SO_q(4)$-isotropic 2-tensor,
deformation of the ordinary Euclidean metric, and
``Killing form'' of $U_qso(4)$.

\medskip
The  commutation relations  and $\star$-conjugation relations are preserved by
the (left) coactions of both $SO_q(4)=SU_q(2)\otimes SU_q(2)'/\b{Z}_2$
and of the extension
$\widetilde{SO_q(4)}:=SO_q(4)\!\times\! GL^+(1)=\b{H}_q\!\times\!\b{H}_q'/GL(1)$
(the quantum group of rotations and scale transformations in 4
dimensions), which take the form
\be
x\to a\, x\,b.                                 \label{SOq4coaction}
\ee
Here $a,b$ are the defining matrices of $SU_q(2),SU_q(2)'$ in the first case and of
$\b{H}_q,\b{H}_q'$
in the second (with entries commuting with each other and with those of $x$),
and matrix product is understood.

\bigskip
A different matrix version (with no interpretation in terms of
$q$-deformed quaternions) of a $SU_q(2)\times SU_q(2)$ covariant
quantum Euclidean space  was proposed in \cite{Maj94}.

\bigskip

 Define 
\bea
&&\alpha'=\sqrt{2}\alpha^{\star} \frac {2}{1\!+\!2|x|^2}e^{ia},\qquad\qquad
\alpha'{}^{\star} =\sqrt{2}\alpha\frac {2}{1\!+\!2|x|^2}e^{-ia},\nn[8pt]
&&\beta'=\sqrt{2}\gamma^{\star} \frac {2}{1\!+\!2|x|^2}e^{ib},\qquad\qquad
\beta'{}^{\star} =\sqrt{2}\gamma\frac
{2}{1\!+\!2|x|^2}e^{-ib},\qquad\label{redef} \\[8pt] && z= \frac
{1\!-\!2|x|^2}{1\!+\!2|x|^2}\nonumber 
\eea
where $\alpha,\beta,...$ fulfill (\ref{explqquatcomrel}) 
and $e^{ia},e^{ib}\in U(1)$ are possible
phase factors. 
Then $\alpha',\beta',z$ fulfill the defining relation (1) of the
$C^{\star}$-algebra considered in Ref. \cite{DabLanMas01}
(where these elements are respectively denoted as $\alpha,\beta,z$),
in particular
\be
\alpha'\alpha'{}^{\star}+\beta'\beta'{}^{\star}+z^2=\1,
\ee
which shows that the noncommutative manifold is a deformation
$S_q^4$ of the 4-sphere.
The invertible function $z(|x|)$ spans
$[-1,1[$, i.e. all the spectrum of $z$ except the eigenvalue $z=1$,
as $|x|$ spans all its spectrum $[0,\infty[$.

The redefinitions
(\ref{redef}) have exactly the form of a stereographic
projection of $\b{R}^4$ on a sphere $S^4$ of unit radius
(recall that $x\cdot x=2|x|^2$): $S^4$ is the sphere centered at the
origin and $\b{R}^4$ the subspace $z=0$ immersing both in a $\b{R}^5$ with
coordinates defined by $X\equiv(Re(\alpha'),Im(\alpha'),
Re(\beta'),Im(\beta'),z)$.
In the commutative theory the  point $X=(0,0,0,0,1)$ of
$S^4$ is the point at infinity of $\b{R}^4$, therefore going
from $\b{R}^4$ to $S^4$ amounts to  compactifying $\b{R}^4$  to $S^4$. 
We can thus regard
the transition from our algebra to the one  considered in Ref.
\cite{DabLanMas01}
as a compactification of $\b{R}_q^4$ into their $S^4_q$.

\section{Other preliminaries}

The $SO_q(4)$-covariant {\bf differential calculus $(d,\Omega^*)$ on
$\b{R}_q^4\sim\b{H}_q$}  \cite{CarSchWat91} is obtained imposing
covariant homogeneous bilinear commutation relations
 (\ref{xxirel}) between the $x^a$ and the differentials
$\xi^a:=dx^a$.
 Partial derivatives are introduced through the
decomposition
$d=\xi^a\partial_a=\xi^{\alpha\alpha'}\partial_{\alpha\alpha'}$. All
other commutation relations are derived by consistency. The complete
list is given by \bea
&& {\sf P}_A{}^{hi}_{jk}x^jx^k=0, \qquad\quad \Leftrightarrow
\qquad\quad x^{\alpha\alpha'}x^{\beta\beta'}=
\hat R^{\alpha\beta}_{\gamma\delta}
\hat R^{-1}{}^{\alpha'\beta'}_{\gamma'\delta'}
x^{\gamma\gamma'}x^{\delta\delta'}\label{xxrel}\\
&& x^h\xi^i=q\hat {\sf R}^{hi}_{jk}\xi^jx^k\qquad\quad \Leftrightarrow
\qquad\quad x^{\alpha\alpha'}\xi^{\beta\beta'}=
\hat R^{\alpha\beta}_{\gamma\delta}
\hat R^{\alpha'\beta'}_{\gamma'\delta'}
\xi^{\gamma\gamma'}x^{\delta\delta'},\label{xxirel}\\
&& ({\sf P}_s+{\sf P}_t)^{ij}_{hk}\xi^h\xi^k=0\qquad \:\Leftrightarrow
\qquad \:\P_s{}^{\alpha\beta}_{\gamma\delta}
\P_s{}^{\alpha'\beta'}_{\gamma'\delta'}\xi^{\gamma\gamma'}
\!\xi^{\delta\delta'}\!=\!0\!=\!(\xi\epsilon\xi^T)^{\gamma\delta}
\epsilon_{\gamma\delta},\label{xixirel}\\
&& {\sf P}_A{}^{ij}_{hk}\partial_j\partial_i=0\qquad\qquad \Leftrightarrow
\qquad \partial_{\alpha\alpha'}\partial_{\beta\beta'}=
\hat R^{\delta\gamma}_{\beta\alpha}
\hat R^{-1}{}^{\delta'\gamma'}_{\beta'\alpha'}
\partial_{\gamma\gamma'}\partial_{\delta\delta'}, \label{ddrel}\\
&& \partial_i x^j = \delta^j_i+q\hat {\sf R}^{jh}_{ik}
x^k\partial_h\quad \: \Leftrightarrow \quad \: \partial_{\alpha\alpha'}
x^{\beta\beta'}\!\!=\!\delta^{\beta}_{\alpha} \delta^{\beta'}_{\alpha'}\!+\! \hat
R^{\beta\delta}_{\alpha\gamma} \hat R^{\beta'\delta'}_{\alpha'\gamma'}
x^{\gamma\gamma'} \!\partial_{\delta\delta'},  \label{dxrel}\\
&& \partial^h\xi^i=q^{-1}\hat {\sf R}^{hi}_{jk}\xi^j
\partial^k\qquad \Leftrightarrow \qquad \partial_{\alpha\alpha'}
\xi^{\beta\beta'}= \hat R^{-1}{}^{\beta\delta}_{\alpha\gamma}
\hat R^{-1}{}^{\beta'\delta'}_{\alpha'\gamma'}\xi^{\gamma\gamma'}
\partial_{\delta\delta'}.\label{dxirel}
\eea
The Laplacian
$\Box\equiv\partial\cdot\partial:=\partial_kg^{hk}\partial_h$ is
$SO_q(4)$-invariant and commutes the $\partial_i$. In ${\cal H}$
there exists a special invertible element $\Lambda$ such that
\[
\Lambda x^i=q^{-1}x^i\Lambda,\qquad\quad
\Lambda\partial^i=q\partial^i\Lambda, \qquad\quad \Lambda
\xi^i=\xi^i\Lambda.\label{Lambdaprop}
\]

{\bf Definitions}:
\begin{itemize}

\item $\bigwedge^*\equiv$ $\natural$-graded algebra generated by the
$\xi^i$, where grading $\natural \equiv$degree in $\xi^i$; any
component $\bigwedge^p$ with $\natural =p$ carries an
irreducible representation of $U_qso(4)$ and has the same dimension
as in the $q=1$ case.

\item ${\cal DC}^*\equiv $ $\natural$-graded algebra generated by
$x^i,\xi^i,\partial_i$. Elements of ${\cal DC}^p$ are
differential-operator-valued $p$-forms.

\item $\Omega^*\equiv$ $\natural$-graded subalgebra generated by the
$\xi^i,x^i$. By definition $\Omega^0=\A$ itself, and both $\Omega^*$
and $\Omega^p$ are $\A$-bimodules. Also, we shall denote $\Omega^*$
enlarged with $\Lambda^{\pm 1}$ as $\tilde\Omega^*$, and the
subalgebra generated by $T^{\alpha\alpha'}:=x^{\alpha\alpha'}/|x|$,
$dT^{\alpha\alpha'}$ as $\Omega_S^*$
(the latter is 4-dim! See below).

\item ${\cal H}\equiv$subalgebra generated by the $x^i,
\partial_i$. By definition, ${\cal DC}^0={\cal H}$, and
both ${\cal DC}^*$ and ${\cal DC}^p$ are ${\cal H}$-bimodules.

\end{itemize}

The special $\widetilde{SO_q(4)}$-invariant 1-form
$$ 
\theta:=\frac 1{1-q^{-2}}|x|^{-2}\,d|x|^2=
\frac{q^{-2}}{q^2-1}\xi^{\alpha\alpha'}\frac{x^{\beta\beta'}}{|x|^2}
\epsilon_{\alpha\beta}\epsilon_{\alpha'\beta'}
$$
plays the role of "Dirac Operator" \cite{Con94} of the differential
calculus,
\[
d\omega_p=[-\theta,\omega_p\}
\equiv -\theta\omega_p+(-)^p\omega_p \theta,\qquad\qquad \omega_p\in\Omega^p ,
\label{thetacommu}
\]
$\theta$ is closed:
\be
d\theta=0,\qquad\theta^2=0.               \label{theta^2=0}
\ee
Applying $d$ to (\ref{blu})  we find
\be
x\bar\xi+\xi\bar x=(q^2\!-\!1)\theta|x|^2I_2,\qquad\qquad
\bar x\xi+\bar\xi x=(q^2\!-\!1)\theta|x|^2 I_2.   \label{xblu}
\ee
Relation (\ref{xxirel}) implies $|x|^2\xi^i=q^2\xi^i|x|^2$, which we generalize
as usual to
\be
|x|^{\pm 1}\xi^i=q^{\pm 1}\xi^i|x|^{\pm 1},  \qquad\Rightarrow\qquad
|x|^{\pm 1}\,\theta=q^{\pm 1}\,\theta\,|x|^{\pm 1}.            \label{theta|x|rel}
\ee
\bigskip
However, $d(f^{\star})\neq (df)^{\star}$, and moreover there is no
$\star$-structure $\star:\Omega^*\to\Omega^*$, but only a
$\star$-structure
$$
\star:{\cal DC}^*\to {\cal DC}^*
$$ \cite{OgiZum92},
with a rather nonlinear character (the latter  has been recently
\cite{Fio04} recast in a much more suggestive form).

\bigskip

The {\bf Hodge map} \cite{Fio94,Fio04JPA} is a $SO_q(4)$-covariant,
$\A$-bilinear map $*:\tilde\Omega^p\to\tilde\Omega^{4-p}$ such that
$*^2= \id$, defined by
\[ 
{}^*(\xi^{i_1}...\xi^{i_p})=q^{-4(p-2)} c_p\,\xi^{i_{p+1}}...\xi^{i_4}
\varepsilon_{i_4...i_{p+1}}{}^{i_1...i_p}\Lambda^{2p-4}, 
\] 
where in our normalization the 
$\varepsilon^{hijk}\equiv$ $q$-epsilon tensor is given by
$$ 
\begin{array}{|c|c|c|c|} 
\hline 
\varepsilon^{-2-112}=q^{-2} & \varepsilon^{-21-12}=-q^{-2} 
&\varepsilon^{-2-121}=-q^{-1} &\varepsilon^{-212-1}=q^{-1} \\  \hline 
\varepsilon^{-22-11}=1 & \varepsilon^{-221-1}=-1 & 
\varepsilon^{-1-212}=-q^{-1} &  \varepsilon^{-11-22}=1 \\  
\hline 
\varepsilon^{-1-221}=1 & \varepsilon^{-12-21}=-1 & 
\varepsilon^{-121-2}=q &  \varepsilon^{-112-2}=-1 \\  
\hline 
\varepsilon^{1-1-22}=-1 & \varepsilon^{1-2-12}=q^{-1} & 
\varepsilon^{1-12-2}=q &  \varepsilon^{12-1-2}=-q \\  
\hline 
\varepsilon^{12-2-1}=1 & \varepsilon^{1-22-1}=-1 & 
\varepsilon^{2-2-11}=-1 &  \varepsilon^{2-1-21}=q \\  
\hline 
\varepsilon^{21-2-1}=-q & \varepsilon^{2-21-1}=1 & 
\varepsilon^{2-11-2}=-q^2 &  \varepsilon^{21-1-2}=q^2 \\  
\hline 
\varepsilon^{-11-11}=k& \varepsilon^{1-11-1}=-k& 
\multicolumn{2}{|c|}{\varepsilon^{ijkl}=0~~~\mbox{otherwise}.} \\ 
\hline 
\end{array} 
$$ 
and $c_p$ are suitable
normalization factors \cite{Fio04JPA}. Actually this extends to a
${\cal H}$-bilinear map $*:{\cal DC}^p\to{\cal DC}^{4-p}$ with the same
features. For $p=2$ the powers of $\Lambda$ disappear and one even gets a map
$*:\Omega^2\to\Omega^2$ defined by 
\be 
{}^*\xi^i\xi^j:=\frac 1{[2]_q}\xi^h\xi^k \varepsilon_{kh}{}^{ij}
=\left(\PH_a-\PH_{a'}\right)^{ij}_{hk}\xi^h\xi^k,
\label{defHodge2x4} 
\ee
where $\PH_a,\PH_{a'}$ were defined in (\ref{defproj})
and $[2]_q=q\!+\!q^{-1}$; the second 
equality can be proved by a direct computation.
$\Omega^2$ (resp. ${\cal DC}^2$) splits into
the direct sum of $\A$- (resp. ${\cal H}$-) bimodules
\[
\Omega^2=\check\Omega^2\oplus \check\Omega^{2}{}' \qquad\quad
\mbox{(resp. }{\cal
DC}^2=\check{\cal DC}^2\oplus \check{\cal DC}^{2}{}'\mbox{)}
\]
of the  eigenspaces of $*$ with eigenvalues $1,-1$ respectively,
whose elements are ``self-dual and anti-self-dual 2-forms''.
$\check\Omega^2$ (resp. $\check{\cal DC}^2$) is generated by the
self-dual exterior forms $(\xi\bar\xi)^{\alpha\beta}$,
or equivalently by the ones
\be f^{\alpha\beta}:=
(\xi\bar\xi\epsilon)^{\alpha\beta}                      \label{deff} \ee
through (left or right) multiplication  by elements of $\A$ (resp.
${\cal H}$). $f^{\alpha\beta}$ span a (3,1) corepresentation
space of $SU_q(2)\times SU_q(2)'$.

One can find 1-form-valued matrices $a$ such that \be
d\,a^{\alpha\beta}=f^{\alpha\beta}; \ee $a$ is uniquely determined
to be 
\be 
a^{\alpha\beta} ={\cal P}_s{}^{\alpha\beta}_{\gamma\delta}
(\xi\epsilon x^T)^{\gamma\delta}, \label{aexplicit} 
\ee
if we require $a^{\alpha\beta}$ to transform as $f^{\alpha\beta}$, i.e. in
the (3,1) dimensional corepresentation of $SU_q(2)\times SU_q(2)'$,
whereas will be defined up to $d$-exact terms of the form
$$
\tilde a= a+\1_2\,dM(|x|^2)
$$  if we  just require $\tilde a^{\alpha\beta}$ to be in the
$(3,1)\oplus (1,1)$ reducible representation. In particular, the
1-form valued matrix
\be
\hat a:=-\xi\bar x,
\ee
as well as the one $(dT)\overline{T}$ (see section
\ref{Wordiffcalc}), belong to the latter,
therefore are invariant  under
the right  coaction of  $SU_q(2)$.
In
the $q=1$ limit (\ref{aexplicit}) becomes
\[
a^{\alpha\beta}=\Big(\xi\epsilon x^T\Big)^{(\alpha\beta)}
=-\left\{Im(\xi\,\bar x\epsilon)\right\}^{\alpha\beta}.
\]

Similarly, antiself-dual $\check \Omega^2{}'$, $\check{\cal
DC}^2{}'$ are generated by $(\bar\xi\xi)^{\alpha'\beta'}$,
or equivalently by
\be
f'{}^{\alpha'\beta'}:=(\bar\xi\xi\epsilon)^{\alpha'\beta'}, \label{deff'}
\ee
and one can find
1-forms $a'{}^{\alpha'\beta'}$ such that
$d\,a'{}^{\alpha'\beta'}=f'{}^{\alpha'\beta'}$, etc.

\bigskip
{\bf Integration over $\b{R}_q^4$} \cite{Ste96,Fio93,fiothesis}
can be introduced by the decompositon
$$
\int_{\b{R}_q^4}d^4x=\int\limits_0^{\infty}d|x| \: \int_{|x|\cdot
S_q^3}d^3T
$$
Integration over the radial coordinate
has to fulfill the scaling property
$\int\limits_0^{\infty}d|x|\,g(|x|)=\int\limits_0^{\infty}d(q|x|)\,g(q|x|)$.
Integration over the quantum sphere $S_q^3$
is determined up to normalization by the requirement of
$SO_q(4)$-invariance. The algebra of functions on  the quantum sphere
$S_q^3$ is generated by the
$T^{\alpha\beta}:=x^{\alpha\beta}/|x|$.

This integration over $\b{R}_q^4$ fulfills all the main properties
of Riemann integration over $\b{R}^4$, including Stokes' theorem,
except the cyclic property, which is $q$-deformed.

\section{Connection with the bicovariant differential calculi on $GL_q(2)$ and
$SU_q(2)$} \label{Wordiffcalc}

We start by recalling that an alternative calculus $(\hat\Omega^*,\hat
d)$ on $\b{R}_q^4$  is obtained by replacing $\hat R\leftrightarrow \hat
R^{-1}$, $q\leftrightarrow q^{-1}$ in relations  (\ref{xxirel}),
$$
x^h\hat\xi^i=q^{-1}\hat {\sf R}^{-1}{}^{hi}_{jk}\hat\xi^jx^k\qquad\quad \Leftrightarrow
\qquad\quad x^{\alpha\alpha'}\hat\xi^{\beta\beta'}=
\hat R^{-1}{}^{\alpha\beta}_{\gamma\delta}
\hat R^{-1}{}^{\alpha'\beta'}_{\gamma'\delta'}
\hat\xi^{\gamma\gamma'}x^{\delta\delta'},
\eqno{(\widehat{\ref{xxirel}})} $$
and in the following ones [(\ref{xxrel})
is invariant under these replacements]. As just done, we shall
add a $\hat {}$  to label these formulae and the corresponding
objects after the replacements.

\bigskip
We first show that
the two differential calculi on $\b{R}_q^4$ coincide with the two 
bicovariant differential calculi on $M_q(2),GL_q(2)$ \cite{S2,SchWatZum92}.
We recall that a differential calculus is completely detemined by the Leibniz
rule and nilpotency for the exterior derivative and by the commutation
relations between the generators of the algebra and their differentials.
For our calculus $(\Omega^*,d)$ the latter read (\ref{xxirel}), whereas for the calculus
on $M_q(2),GL_q(2)$ they are (13)$_1$ in \cite{SchWatZum92}.
Now it is straightforward to check that indeed relation (\ref{xxirel}), in the matrix
formulation at the right, amounts to relation (13)$_1$ in
\cite{SchWatZum92}, provided we identify $x\to A$ and recall that $\hat R:= P
R$ ($P$ deonting the permutation matrix), $\hat R^T=\hat R$. To complete the
`dictionary' we add that our $T,\theta, \bar\xi x$ have to be identified with
$T,(q^{-1}-q)^{-1}\xi,-\Omega$ of \cite{SchWatZum92}.

\bigskip
We now verify that, restricting as in \cite{SchWatZum92} either calculus to the
subalgebra generated by the $T^{\alpha\alpha'}=x^{\alpha\alpha'}/|x|$, one
obtains differential calculi $(\Omega_S^*,d)$, $(\hat\Omega_S^*,\hat d)$ on
$SL_q(2)$, which coincide with Woronowicz $4D\mp$ bicovariant differential
calculi \cite{Wor89,PodWor90}. [For real $q\neq 0,1$ the latter are also real,
i.e compatible with $(da)^{\star}=d(a^{\star})$ and
the $\star$-structure (\ref{starSUq2}) of $SU_q(2)$.]

Introduce the 1-form valued matrix $\omega:=\xi\bar x/|x|^2$. Using
(\ref{xxrel}), (\ref{xxirel}), (\ref{theta|x|rel}), $\hat R^T=\hat R$ and $$
\epsilon_{\alpha\lambda}\hat R^{\pm 1}{}_{\beta\gamma}^{\lambda\mu}
=q^{\pm 1}\hat R^{\mp 1}{}_{\alpha\beta}^{\mu\lambda}
\epsilon_{\lambda\gamma}
$$
[which is a consequence of (\ref{qepsilon})] it is easy to show that
 \be
T^{\alpha\alpha'}\omega^{\beta\beta'}
=q^{-1}\hat R^{\alpha\beta}_{\lambda\delta}
\hat R^{\mu\delta}_{\gamma\beta'}\,
\omega^{\lambda\mu}\,T^{\gamma\alpha'}.          \label{Txixrel}
\ee
On the other hand, by a straightforward computation one finds
$$
 dT^{\alpha\alpha'}=
q^{-1}\xi^{\alpha\alpha'}\frac 1{|x|}+(q^{-1}\!-\!1)\theta T^{\alpha\alpha'},
$$
whence
\be
(dT)\overline{T}  =q^{-1}\omega+(q^{-1}\!-\!1)\theta I_2.
\label{RLinv1forms}
\ee
This 1-form-valued `Maurer-Cartan' $2\times 2$ matrix and the one
$(d\overline{T})T$
are by (\ref{SOq4coaction}) manifestly invariant under respectively the right
and left coaction of  $SU_q(2)$, or equivalently under the $SU_q(2)'$ and  the
 $SU_q(2)$ part of $SO_q(4)$ coaction. Setting
$Q:=-\epsilon^{-1}\epsilon^T$  one finds
\be
\mbox{tr}[Q(dT)\overline{T}]=\mbox{tr}[Q^{-1}(d\bar
T)T]=(q\!-\!1)(q\!-\!q^{-2}) \theta;             \label{thetaTT}
\ee
only in the $q\to 1$ limit these traces vanish. That's why for generic
$q\neq 1$ the four matrix elements of either $(dT)\overline{T}$ or
$(d\overline{T})T$ are independent (4-dimensional calculus) and make up
alternative bases for  $\Omega^*$.
Moreover, we see that for $q\neq 1$ the `Dirac operator' $\theta$ can be
expressed purely in terms of the matrix elements of $dT$ and $T$,
in other words the restriction $(d,\Omega^*_S)$ of the above calculus
to $C\big(SU_q(2)\big)$ is well defined and 4-dimensional.
From (\ref{RLinv1forms}) one sees that the
matrix elements $\omega^{\beta\beta'}$ make up an alternative basis of
$\Omega^*_S$; their
commutation relations (\ref{Txixrel}) with the $T^{\alpha\alpha'}$
completely specify the first calculus.
Similarly, setting for the other calculus $\hat\omega:=\hat\xi\bar x/|x|^2$ we
find
$$
T^{\alpha\alpha'}\hat\omega^{\beta\beta'}
=q\hat R^{-1}{}^{\alpha\beta}_{\lambda\delta}
\hat R^{-1}{}^{\mu\delta}_{\gamma\beta'}\,
\hat\omega^{\lambda\mu}\,T^{\gamma\alpha'},   \eqno{(\widehat{\ref{Txixrel}})}
$$
$$
(\hat dT)\overline{T}  =q\hat \omega+(q\!-\!1)\hat\theta I_2,
                                   \eqno{(\widehat{\ref{RLinv1forms}})}
$$
and
$$
\mbox{tr}[Q(\hat dT)\overline{T}]=\mbox{tr}[Q^{-1}(\hat d\bar
T)T]=(q\!-\!1)(q\!-\!q^{-2}) \hat \theta.     \eqno{(\widehat{\ref{thetaTT}})}
$$

Let us compare now our results with
Woronowicz $4D+$  bicovariant differential calculus on
$C\big(SU_q(2)\big)$ \cite{Wor89,PodWor90}. We describe the latter in the $R$-matrix
formalism, as done in Ref. \cite{CarSchWatWei91}, where the matrix $T$ was
denoted as $M$. Comparing formula (5.8) of the latter with our
($\widehat{\ref{thetacommu}}$) leads to identify our bi-invariant 1-form `Dirac
operator' $\hat \theta$ with their $ -X/{\cal N}$. This is consistent as we
then find that our ($\widehat{\ref{thetaTT}}$) coincides with their (5.26)
(with $N=2$).  Formula (5.23) of \cite{CarSchWatWei91}  [$\kappa$ denotes the
antipode, so $\kappa(M)$ is our $T^{-1}=\overline{T}$] leads to
identify our right invariant 1-form valued matrix $(\hat dT)\overline{T}$ with
their $\tilde\theta$; further comparison of formula (5.25) of
\cite{CarSchWatWei91} with our ($\widehat{\ref{RLinv1forms}}$)
leads to identify our $\omega^{\alpha\alpha'}$ with their
$\theta^{\alpha}_{\alpha'}(1-q^2)/{\cal N}q^3$ (but they use latin letters
instead of greek ones to label matrix rows and columns). This is consistent
because the commutation relations ($\widehat{\ref{Txixrel}}$)
coincide with the commutation relations for the $\theta^{\alpha}_{\alpha'}$
which one obtains after little work from their formulae (4.14), (3.16), (3.20)
and the  $\star$-conjugates of the latter. Therefore the differential calculus
$(\hat\Omega_S^*,\hat d)$ coincides with Woronowicz $4D+$
bicovariant one. Similarly one shows that the differential calculus
$(\Omega_S^*,d)$ coincides with  Woronowicz $4D-$
bicovariant one.

\bigskip
We end by noting that the above identifications and our results about the
Hodge  map give as a bonus a well-defined Hodge operator $*$ and (anti)selfdual
2-forms on $M_q(2),GL_q(2),SL_q(2),SU_q(2)$. On $M_q(2),GL_q(2)$ (anti)selfdual
2-forms are respectively the $(\xi\bar\xi)^{\alpha\beta}$,
$(\bar\xi\xi)^{\alpha'\beta'}$, whereas on $SL_q(2),SU_q(2)$  are
respectively obtained dividing
$(\xi\bar\xi)^{\alpha\beta}$,
$(\bar\xi\xi)^{\alpha'\beta'}$ by $|x|^2$ and
expressing the results in term of $T, dT$ only:
\bea
&&v^{\alpha\beta}:=(\xi\bar\xi)^{\alpha\beta}\frac {q^{-1}}{|x|^2}=
\left[q^2T\theta\overline{T}\theta+\theta T\theta\overline{T}\right]^{\alpha\beta}=
\left[q^2(dT)(d\overline{T})+(q^2-1)(dT)\theta\overline{T}\right]^{\alpha\beta}
\qquad\qquad\\
&&v^{\alpha'\beta'}:=(\bar\xi\xi)^{\alpha'\beta'}\frac {q^{-1}}{|x|^2}=
\left[q^2\overline{T}\theta T\theta+\theta \overline{T}\theta T\right]^{\alpha'\beta'}=
\left[q^2(d \overline{T})(dT)+(q^2-1)(d\overline{T})
\theta T\right]^{\alpha'\beta'}\qquad\qquad
\eea
fulfill
\be
{}^*v^{\alpha\beta}=v^{\alpha\beta},\qquad\qquad {}^*v^{\alpha'\beta'}=-v^{\alpha'\beta'}
\ee
[and similarly for the other calculus $(\hat\Omega_S^*,\hat d)$].

\section{Formulations of noncommutative gauge theories}

We recall some minimal common elements in the formulations of $U(n)$ gauge theories
on commutative as well as noncommutative spaces
\cite{Con94,Lan97,FigGraVar01,Mad99}. In $U(n)$ gauge
theory the gauge transformations $U$ are unitary $\A$-valued ($\A$
being the algebra of functions on the noncommutative manifold)
$n\times n$  matrices, $U\!\in\! M_n(\A)\equiv
M_n(\b{C})\otimes_{\b{C}}\A$. The gauge potential $A\equiv (
A^{\dot\alpha}_{\dot\beta})$ is a antihermitean 1-form-valued
$n\times n$ matrix, $A\in M_n(\Omega^1(\A))$. The definition of the
field strength $F\in M_n(\Omega^2(\A))$ associated to $A$ is as
usual $F:=dA+AA$. At the right-hand side the product $AA$ has to be
understood both as a (row by column) matrix product and as a wedge
product. Even for $n=1$, $AA\neq 0$, contrary to the commutative
case. The Bianchi identity $DF:=dF+[A,F]=0$ is automatically
satisfied and the Yang-Mills equation reads as usual $D{}^*F=0$.
Because of the Bianchi identity, the latter is automatically
satisfied by any solution of the (anti)self-duality equations \be
{}^*F=\pm F. \ee

The Bianchi identity, the  Yang-Mills equation, the
(anti)self-duality equations,  the
flatness condition $F=0$ are preserved by gauge transformations
$$
A^U=U^{-1}(AU +dU), \qquad\Rightarrow \qquad F^U= U^{-1}F U.
$$
As usual, $A=U^{-1}dU$ implies $F=0$. Up to normalization factors,
the gauge invariant `action' $S$ and `Pontryagin index'
(or `second Chern number') $\c{Q}$
are defined by 
\be
S = \mbox{Tr}(F\:{}^*\!F), \qquad\qquad\qquad  \c{Q} = \mbox{Tr}(FF)
              \label{actionfun} 
\ee 
where Tr stands for  a
positive-definite trace combining the $n\times n$-matrix trace with
the integral over the noncommutative manifold (as such, Tr has to
fulfill the cyclic property). If integration $\int$ fulfills itself
the cyclic property then this is obtained by simply choosing $
\mbox{Tr}=\int  \mbox{tr}$, where $\mbox{tr}$ stands for the
ordinary matrix trace. $S$ is automatically nonnegative.

In commutative geometry the socalled Serre-Swan theorem \cite{Ser62,Con95}
states that vector bundles over a compact manifold coincide
with finitely generated projective modules
${\cal E}$ over $\A$. The gauge connection $A$ of a gauge group
(fiber bundle) acting on a vector bundle is expressed in terms
of the projector ${\cal P}$ characterizing the projective module. Therefore
these projectors can be used to  completely determine the connections.
In Connes' standard approach \cite{Con94} to noncommutative geometry
the finitely generated projective modules are the primary objects to
define and develop the gauge theory. The topological
properties of the connections can be classified in terms of topological
invariants (Chern numbers), and the latter can be computed
directly in terms of characters of  ${\cal P}$ (Chern-Connes
characters), in particular $\c{Q}$ can be computed in terms of 
the second Chern-Connes character, when
Connes' formulation of noncommutative geometry applies.

\bigskip

In the present $\A\equiv C(\b{R}_q^4)=C(\b{H}_q)$ case there are {\bf 2 main
problems} preventing the application of this formulation of gauge theories:

\begin{enumerate}

\item Integration over $\b{R}_q^4$ fulfills a {\it deformed} cyclic property
\cite{Ste96}.

\item $d(f^{\star})\neq (df)^{\star}$, and there is no
$\star$-structure $\star:\Omega^*\to\Omega^*$, but only a
$\star$-structure  $\star:{\cal DC}^*\to {\cal DC}^*$ \cite{OgiZum92},
with a nonlinear character.

\end{enumerate}

A solution to both problems might be obtained
\begin{enumerate}

\item allowing for ${\cal DC}^1$-valued $A$ ($\Rightarrow$
${\cal DC}^2$-valued $F$'s), and/or

\item realizing Tr$(\cdot)$
by in the form $\mbox{Tr}(\cdot)$:=$\int\mbox{tr}(W\cdot)$, with $W$ some
suitable positive definite ${\cal H}$-valued (i.e.
pseudo-differential-operator-valued) $n\times n$ matrix (this
implies a change in the hermitean conjugation of differential
operators), or even a more general form.

\end{enumerate}

This hope is based on our results \cite{Fio04}.

\section{The (anti)instanton solution}

 We first recall the {\it commutative ($q\!=\!1$)
solution} of the self-duality eq. ${}^*F=F$: the
instanton solution of \cite{BelPolSchTyu75} in t' Hooft
\cite{tHo76} and in ADHM \cite{AtiDriHitMan78}  quaternion notation
(see \cite{Ati79} for an introduction) reads: \bea A &=&  dx^i\, \sigma^a\,
\underbrace{\eta^a_{ij}x^j \frac 1{\rho^2+r^2/2}}_{A^a_i},\nn
&=&-Im\left\{\xi\,\frac{\bar x}{|x|^2}\right\} \frac
1{1+{\rho^2}\frac 1{|x|^2}}\nn  &=& -(dT)\overline{T} \frac
1{1+{\rho^2}\frac 1{|x|^2}}\label{inst}\\
 F &=& \xi\bar\xi\,\rho^2\frac
1{(\rho^2+|x|^2)^2}, \label{Finst} \eea where $r^2:=x\cdot
x=2|x|^2$, $\eta^a_{ij}$ are the so-called 't Hooft $\eta$-symbols and $\rho$
is the size of the instanton (here centered at the origin).
The third equality is based on the identity
$$
\xi\,\frac{\bar x}{|x|^2}=(dT)\overline{T}+I_2\frac {d|x|^2}{2|x|^2}
$$
and the observation that the first and second term at the rhs are
respectively antihermitean and hermitean,  i.e. the imaginary and
the real part of the quaternion at the lhs.

\bigskip
{\it Noncommutative (i.e. $q\!\neq \!1$) solutions} of ${}^*F=F$. Looking for
$A$ directly in the form  $A=\xi\bar x\,l/|x|^2+\theta\, I_2\,n$,
where $l,n$ are functions of $x$ only through $|x|$,
one finds a family of solutions
parametrized by  $\rho^2$ (a nonnegative constant, or more generally a further
generator of the algebra) and by the function $l$ itself.
The freedom in the choice of $l$ should disappear upon
imposing the proper (and still missing!) antihermiticity
condition on $A$, as it occurs in the $q=1$ case. For the moment, out of this
large family we just pick one which has the right $q\to 1$ limit
and closely resembles the undeformed solutions (\ref{inst}-\ref{Finst}):
\be
\ba{l} A = - (dT)\overline{T} \,\frac 1{1+{\rho^2}\frac
1{|x|^2}},\\[8pt] F = q^{-1}\xi\bar\xi \frac 1{|x|^2+\rho^2}
\rho^2\frac 1{q^2|x|^2+\rho^2}.\ea\label{qinst}
\ee
Of course we have to extend the algebras so that they contain the rational
functions at the rhs.
The matrix elements $A^{\alpha\beta}$  span a
$(3,1)\oplus (1,1)$ dimensional corepresentation of
$SU_q(2)\times SU_q(2)'$,
suggesting as the `fiber' of
the gauge group  in the complete theory
a (possibly deformed) $U(2)$ [instead of a $SU(2)$].

By the scaling and translation invariance of integration over
$\b{R}_q^4$,  if we could find a `good' pseudodifferential operator
$W$  to define gauge invariant ``action'' and ``topological charge''
by $$ \c{Q}:=\int_{\b{R}_q^4}\mbox{tr}(WF\,F)=
\int_{\b{R}_q^4}\mbox{tr}(W F\,{}^*F)=S $$ the latter would, as in
the commutative case, equal a constant independent of $\rho,y$
(which by the choice of the normalization of the integral we can
make 1).

In the $q=1$ case multi-instanton solution are explicitly written
down in the socalled `singular gauge'. Note that as in the $q=1$
case $T=x/|x|$ is unitary and singular at $x=0$. So it can play the
role of a `singular gauge transformation'. In fact $A$ can be
obtained through the singular gauge transformation $A=T(\hat
A\overline{T}+d\overline{T})$ from the singular gauge potential 
\bea
\hat A &=& \overline{T}dT\frac 1{1+|x|^2\frac 1{\rho^2}} =-\frac
1{1+|x|^2\frac 1{q^2\rho^2}}(d\overline{T})T\label{hatA0}\\
&=&-\frac 1{1+|x|^2\frac 1{q^2\rho^2}} \left[q^{-1}\bar\xi\frac x{|x|^2}-
\frac {q^{-3}I_2}{1\!+\!q}\left(
\xi^{\alpha\alpha'}\frac{x^{\beta\beta'}}{|x|^2}
\epsilon_{\alpha\beta}\epsilon_{\alpha'\beta'}\right)\right]. \label{hatA}
\eea
$\hat A$ can be expressed also in
the form \[
\hat A= \phi^{-1}\hat{\cal D}\phi, \qquad\qquad \phi:=1+q^2\rho^2\frac
1{|x|^2},
\]
where $\hat{\cal D}$ is the first-order-differential-operator-valued
$2\times 2$ matrix obtained from the square bracket in (\ref{hatA}) by the
replacement $x^{\alpha\alpha'}/|x|^2\to q^2\partial^{\alpha\alpha'}$:
\be
\hat{\cal D}:=q\bar\xi\partial-\frac {q^{-1}I_2}{q\!+\!1}d \label{Doper}
\ee
(for simplicity we are here assuming that $\rho^2$ commutes with
$\xi^{\alpha\alpha'}\partial^{\beta\beta'}$). $\phi$ is harmonic:
$$
\Box\phi=0.
$$
This is the analog of the $q=1$ case, and is useful for the
construction of multi-instanton solutions.

\bigskip

The {\bf anti-instanton solution} is obtained just by converting
unbarred into barred matrices, and conversely, as in the $q=1$ case.
For instance, from (\ref{qinst}) we obtain the anti-instanton
solution in the regular gauge \be \ba{l} A' = - (d\overline{T})T
\,\frac 1{1+{\rho^2}\frac
1{|x|^2}},\\[8pt] F' = q^{-1}\bar\xi\xi \frac 1{|x|^2+\rho^2}
\rho^2\frac 1{q^2|x|^2+\rho^2}.\ea\label{qantiinst}
\ee

\subsection*{Recovering the instanton projective module of Ref.
\cite{DabLanMas01}} \label{compare}

In commutative geometry
the instanton projective module ${\cal E}$ over $\A$ and the associated gauge
connection can be most easily obtained using the quaternion formalism, in the
way described e.g. in Ref. \cite{Ati79}. $\b{H}\sim\b{R}^4$
can be compactified as $P^1(\b{H})\sim S^4$. Let $(w,x)\in\b{H}^2$ be
homogenous coordinates of the latter, and choose $w=I_2$ on the chart
$\b{H}\sim\b{R}^4$. The element  $u\in\b{H}^2$ defined by
\be
u \equiv \left(\ba{c} u_1\\ u_2\ea \right)=
\left(\ba{c} I_2\\ \frac{\rho x}{|x|^2}\ea \right)
\left(1\!+\!\frac{\rho^2}{|x|^2}\right)^{-1/2}   \label{defu}
\ee
fulfills $u^{\dagger}u=I_2\1$, and the $4\times 2$ $\A$-valued  matrix
$u$ has only three independent components.
Therefore the $4\times 4$ $\A$-valued matrix
\be
{\cal P}:=uu^{\dagger}=\left(\ba{lll} I_2 & \frac{\rho \bar x}{|x|^2}\\
\frac{\rho x}{|x|^2} & \frac{\rho^2}{|x|^2}I_2 \ea \right)
\frac 1{1\!+\!\frac{\rho^2}{|x|^2}}                    \label{calP}
\ee
is a self-adjoint three-dimensional projector. It is
the projector associated in the Serre-Swan theorem
correspondence to the gauge connection (\ref{hatA0}), by the formula $\hat
A=u^{\dagger}du$.
The associated projective module ${\cal E}$ is embedded
in the free module $\A^{16}$ seen as $M_4(\A)$, and is obtained from the latter
as ${\cal E}={\cal P}M_4(\A)$.

In the present $q$-deformed setting we immediately check that
the element $u\in\b{H}_q^2$ defined by (\ref{defu}) fulfills
$u^{\dagger}u=I_2\1$ again, so that
the $4\times 2$ $\A$-valued  matrix ${\cal P}$ defined by (\ref{calP})
is again hermitean and idempotent, and has only 3 independent components.
Therefore, it defines the `instanton projective module'
${\cal E}={\cal P}M_4(\A)$ also in the $q$-deformed case.
One can easily verify that ${\cal P}$ reduces to the hermitean idempotent
$e$ of \cite{DabLanMas01} if one chooses the instanton
size as $\rho=1/\sqrt{2}$ and performs the change of generators
(\ref{redef}). Therefore, interpreting the model  \cite{DabLanMas01}
as a compactification to $S_q^4$ of ours, we can use all
the results \cite{DabLanMas01} about the Chern-Connes classes of $e$.

Unfortunately in the  $q$-deformed case it is no more true that
$\hat A=u^{\dagger}du$,
essentially because the $|x|$-dependent global factor multiplying
the matrix at the rhs(\ref{calP}) does not commute with the 1-forms
of the present calculus ($|x|\xi^i=q\xi^i|x|$).

\subsection*{Shifting the `center of the instanton' away from the origin}
This can be done by
the replacement (or `braided coaddition' \cite{Maj95})
$$
x\to x-y,
$$
where the `coordinates of the center' $y^i$ generate a new copy of
$\A$, `braided' with the original one (see below).
Therefore the instanton moduli space must be a noncommutative manifold, with
coordinates $\rho, y^i$! This is similar to what was proposed
in \cite{IvaLecMue04} for the instanton moduli space  on $\b{R}_{\theta}^4$.
This shift also changes the gauge transformation relating $A,\hat A$
as follows
$$
T=\frac x{|x|}\to \frac {x-y}{|x-y|},
$$
namely we must now allow also gauge transformations depending on the
additional noncommutative parameters.

\section{Multi-instanton solutions}

We have found solutions of the self-duality equation corresponding to
$n$ instantons in the ``singular
gauge''  \cite{tHo76,vari} in the form
\be
\hat A= \phi^{-1}
\hat{\cal D}\phi,
\ee
where $\phi$ is the harmonic scalar function
\be
\phi=1+\rho_1^2\frac 1{(x\!-\!y_1)^2}+
\rho_2^2\frac 1{(x\!-\!y_1\!-\!y_2)^2}+...+\rho_n^2\frac 1{(x\!-\!y_1\!-\!...\!-\!y_n)^2}
\ee
as in the commutative case. In the commutative limit
\bea
&& \rho_{\mu}\equiv \mbox{size of the $\mu$-th instanton},\nn
&& v^i_{\mu}:=\sum\limits_{\nu=1}^{\mu}y^i_{\nu}
\equiv \mbox{$i$-th coordinate of the $\mu$-th instanton}.\nonumber
\eea
are constants ($\mu=1,2,...,n$).
In the noncommutative setting the new generators
$\rho_{\mu}^2,y^i_{\nu}$ have to fulfill the following nontrivial
commutation relations: \bea
&&\rho_{\nu}^2\rho_{\mu}^2=q^2\,\rho_{\mu}^2\rho_{\nu}^2\qquad
\qquad \nu<\mu             \nn
&&\rho_{\nu}^2y_{\mu}^i=y_{\mu}^i\rho_{\nu}^2 \cdot \left\{\ba{l}
q^{-2}\:\: \nu<\mu \cr 1, \:\:\nu \ge\mu \ea\right. \nn
&&\rho_{\mu}^2\xi^i=\xi^i\rho_{\mu}^2, \qquad\qquad
\partial_i\rho_{\mu}^2=\rho_{\mu}^2\partial_i.\label{ntcr}\\
&& y_{\mu}^iy_{\nu}^j=q\RH^{ij}_{hk}y_{\nu}^hy_{\mu}^k\qquad\qquad
\nu<\mu, \nn &&{\sf P\!}_A{}^{ij}_{hk}y_{\mu}^hy_{\mu}^k=0. \nonumber
\eea
($\mu,\nu=0,1,...,n$, and we have set $x^i\equiv y_0^i$).

The last relation states that for any fixed $\nu$ the 4 coordinates
$y_{\nu}^i$ generate a copy of $\A$. The last but one relation states that
the various copies of $\A$ are {\it braided} \cite{Maj95} w.r.t.
each other (this is necessary for the $SO_q(4)$ covariance of the
overall algebra).

\medskip
The obvious consequence of the nontrivial
commutation relations (\ref{ntcr}) is that in a complete theory
{\bf the instanton moduli space must be a noncommutative manifold}.

\medskip
Not only for $n=1$, but also for $n=2$ we have been able to go to a gauge
potential $A$ `regular'  in $z_{\mu}^i\!:=\!x^i\!-\!v^i_{\mu}$
by a `singular gauge transformation', which also depends on
$y^i_{\nu}$ (as in the $q=1$ case
\cite{GiaRot77,OliSciCre79,vari}):
\be
A_2=U_2^{-1}\left(\hat A U_2+dU_2\right),\qquad\qquad      
U_2\equiv U_2(z_1,z_2):=\frac{\bar z_1}{|z_1|}\frac{y_2}{|y_2|}
\frac{\bar z_2}{|z_2|}     \label{A_2}
\ee

\subsection*{Acknowledgments}
The author acknowledges 7 days hospitality and support by the Erwin Schroedinger Institute of Vienna within the program ``Gerbes, Groupoids, and Quantum Field Theory'', during which this work was finished. He is grateful to P. Aschieri and H. Steinacker for stimulating discussions.

\end{document}